\documentclass[a4paper, 12pt]{article}

\usepackage{algorithm}
\usepackage{color,cite,soul}
\usepackage{hyperref}
\usepackage{algorithmic}

\usepackage{graphicx}
\usepackage{stix}
\usepackage{amsmath,amsfonts,amstext,amssymb,amsbsy,amsopn,amsthm,eucal}
\usepackage[english]{babel}

\date{}

\newcommand{\re}{{\mathcal Re}}

\renewcommand{\epsilon}{\varepsilon}
\renewcommand{\tilde}{\widetilde}
\renewcommand{\hat}{\widehat}

 \newcommand{\trace}{\mathop{\mathrm{trace}}}
 \newcommand{\diag}{\mathop{\mathrm{diag}}}
 \newcommand{\argmin}{\mathop{\mathrm{arg\,min}}}
 \newcommand{\argmax}{\mathop{\mathrm{arg\,max}}}
\newcommand\C{\ensuremath{\mathbb{C}}}
\newcommand{\Cn}{\ensuremath{\C^{n}}}

\newcommand{\cH}{{\mathcal{H}}}

\newtheorem{example}{Example}[section]

\DeclareMathOperator{\orth}{orth}

\date{}

\title{Damping optimization of parameter dependent mechanical systems by rational interpolation}

\author{Zoran Tomljanovi\'c\thanks{Department of Mathematics,
  University J.J. Strossmayer in Osijek,  Osijek, Croatia,
     {\tt ztomljan@mathos.hr}}
       \and
Christopher Beattie\thanks{Department of Mathematics, Virginia Polytechnic Institute and State University, Blacksburg,  USA, {\tt
  \{beattie, gugercin\}@vt.edu}}
  \and
     Serkan Gugercin\footnotemark[2]
}

\begin{document}

\maketitle

\begin{abstract}
 We consider an optimization problem  related to  semi-active damping of vibrating systems. The main problem is to determine the best damping matrix  able to  minimize influence of the input on  the output of the system. We use a minimization criteria based on the $\mathcal{H}_2$   system norm.
 The objective function is non-convex  and the  associated optimization problem typically requires a  large number of objective function evaluations.  We propose an optimization approach that calculates `interpolatory' reduced order models, allowing for significant acceleration of the optimization process.

In our approach, we use parametric model reduction (PMOR)  based on the  Iterative Rational Krylov Algorithm, which ensures good approximations relative to the  $\mathcal{H}_2$  system norm, aligning well with the underlying damping design objectives. For the parameter sampling that occurs within each  PMOR cycle, we consider approaches with predetermined sampling and approaches using adaptive sampling, and each of these approaches may be combined with three possible strategies for internal reduction.  In order to preserve important system properties, we maintain second-order structure, which through the use of modal coordinates, allows for very efficient implementation.

The methodology proposed here provides a significant acceleration of the optimization process; the gain in
efficiency is illustrated in numerical experiments.
\end{abstract}
 {\bf msc2010:}{ 49J15, 74P10, 70Q05, 41A05}

\section{Introduction}
We consider the following  vibrational system
 described by
\begin{equation} \label{MDK}
\begin{array}{rl}
M\ddot q(t)+C\dot q(t)+Kq(t)&=E w(t),\\
z(t)&=Hq(t)
\end{array}
\end{equation}
where the mass matrix, $M$, and stiffness matrix, $K$,
are real, symmetric positive-definite matrices of order $n$.
The state variables comprise displacement and rotational degrees of freedom and are collected in the coordinate vector $q\in \mathbb{R}^n$.
  The vector $z\in \mathbb{R}^{m_{out}}$ denotes the
observed output comprising the system displacements in which we are interested, which in turn determines the (constant) output port matrix
$H \in \mathbb{R}^{m_{out} \times n}$.
The vector  $w(t)\in \mathbb{R}^{m_{in}}$ represents a dynamic disturbance, the primary excitation . The associated locations of
the primary excitation determine the
primary excitation  matrix, $E\in \mathbb{R}^{n\times m_{in} }$.

The damping matrix is given by $$C=C_{int} + C_{ext},$$ where $C_{int}$ and $C_{ext}$  represent contributions from internal and external damping, respectively.  We are principally concerned with optimizing the placement and geometry of external dampers whose dynamic effects are described through   $C_{ext}$.   The internal damping contribution as encoded in $C_{int}$ is typically
 both small in magnitude and difficult to model in detail.  Often it is modeled
by taking $C_{int}$ to be a small multiple of critical damping (see,  e.g.,  \cite{BennerTomljTruh10, BennerTomljTruh11,  TRUHVES09}):
\begin{equation}\label{C_int}
  C_{int} = 2 \alpha_c M^{1/2}\left(M^{-1/2}KM^{-1/2}\right)^{1/2}M^{1/2}, \qquad \alpha_c\ll 1.
\end{equation}
Other models of internal damping may be considered as well, but whatever model may be chosen (i.e., however $C_{int}$ is determined), it is part of the description of the underlying structure and not accessible to modification in the course of optimization.

We assume that the external damping, the component of damping that is accessible to modification and optimization,
 creates dynamic effects that may be modeled as $C_{ext}=BGB^T$ where
  $G=\diag{(g_1, g_2, \ldots, g_p)}\in \mathbb{R}^{p\times p}$ is a diagonal matrix  and $B \in \mathbb{R}^{n\times p}$
  determines the placement and geometry of the external dampers.
The  entries $\{g_i\}_{i=1}^p$ are usually called gains or viscosities and represent
friction coefficients  of the corresponding dampers. These coefficients are non-negative and may be either constant or
vary over time.  Here, we consider the $g_i$ to have constant values that will be chosen within fixed feasible margins, $g_i \in  [g_i^-,g_i^+]$ for $i=1,\ldots,p$. Usually the number of dampers $p$ is much smaller than the full dimension: $p\ll n$.
More details regarding system stability and model description can be found in \cite{BennerKuerTomljTruh15, Blanchini12}.

 Damping optimization has a long history in the service of structural engineers and these problems continue to be widely investigated from the perspectives of both engineering and mathematics.  A frequent context for damping optimization comes from system stabilization goals, wherein one strategically adds damping to a structure with light internal damping so as to mute resonances or move them away from the frequencies of ambient oscillatory loads or to suppressing damaging effects of external impacts on the structure. There is a vast literature in this field of research, see, e.g., \cite{MullerSchiehlen85, IA80, NR89, Tak97, Kan13}.
 Depending on the particular application, different damping criteria may be appropriate. For example, for stationary systems criteria that involve spectral abscissa are useful (see \cite{FreitLanc:99}), while criteria that involve total average energy was used in, e.g.,\cite{TrTomPuv16, TRUHVES09}.  Structured dimension reduction methods using total average energy criterion were considered in \cite{BennerTomljTruh10, BennerTomljTruh11}. For non-stationary systems one may consider in addition particular external forces that potentially play an important role in system behaviour.  In this case,  criteria that involve average energy amplitude and average displacement amplitude can be considered (see, e.g., \cite{KuzmTomljTruh12, TRUTOMVES2014, Kan13}).  Overviews on different possible damping criteria can be found in \cite{Gaw, VES2011}.

Since optimization of gains $g_1, g_2, \ldots, g_p$ will  be our main interest, we collect these
parameters into a vector $g=(g_1, g_2, \ldots, g_p)$ and write $G(g)$ to express the parametric dependence on gains.
The transfer function matrix for system (\ref{MDK}) is given by
\begin{equation}\label{TF}
F(s;g)=H (s^2M+s(C_{int} +BG(g)B^T)+K)^{-1} E,\quad s\in \mathbb{C} .
\end{equation}
Note that for any $s\in \mathbb{C}$, $F(s;g)$ is an $m_{out}\times m_{in}$ complex matrix with rational functions of $s$ as elements.
We may also rewrite the second-order system representing our vibrating system as a first-order system of differential equations:
\begin{align}
\dot x(t)&=\widehat{A} x(t)+\widehat{E}w(t),  \label{systemLinearized} \\
z(t)&=\widehat{H}x(t)\nonumber,
\end{align}
where
\begin{equation} \label{matrixAHE}
\begin{array}{c}
x=\left[%
\begin{array}{c}
  q \\
  \dot{q} \\
\end{array}%
\right], \qquad\  \widehat{A}(g)=\left[%
\begin{array}{cc}
  0 & I \\
  -M^{-1}K & -M^{-1}(C_{int}+BG(g)B^T) \\
\end{array}%
\right],\\[4mm]
\hspace*{-10mm} \widehat{E}=\left[%
\begin{array}{c}
  0 \\
  M^{-1}E \\
\end{array}%
\right], \quad  \mbox{and} \quad  \widehat{H}=\left[
         \begin{array}{cc}
           H & 0 \\
         \end{array}
       \right],
       \end{array}
\end{equation}
leading to an alternate realization of $F(s;g)$: $\displaystyle F(s;g)=\widehat{H}\left(s\, I-\widehat{A}(g)\right)^{-1}\widehat{E}$.
The main task in our setting will be to determine a damping configuration, as encoded in
$G$ and $B,$ which minimizes the influence of the input disturbance, $w$,  on observed output states, $z$.
%

%
One may consider different criteria to achieve this goal,
some are based on system-theoretic norms (see, e.g., \cite{Gaw, Blanchini12, BennerKuerTomljTruh15}).  We focus on the $\mathcal{H}_2$ norm and define a cost function in the frequency domain using the  transfer function matrix defined above in \eqref{TF}:
\begin{equation}\label{H2critwithTF}
\left\| F( \cdot \,;g)  \right\|_{\cH_2} = \left (\frac{1}{2\pi} \int_{-\infty}^{+\infty} \trace{(F(j\omega\,;g)^*F(j\omega\,;g))}  \,d\omega \right )^{\frac{1}{2}}.
\end{equation}
 For single-input/single-output (SISO) systems, this criterion can be identified with the response energy resulting from an impulsive input:
\begin{equation}
\left\| F( \cdot \,;g)  \right\|_{\cH_2}^2= \int_0^{+\infty} \|z_D(t)\|_2^2 \,dt.
\end{equation}
Here, $z_D$ is the (scalar) output response of the SISO system (\ref{systemLinearized}) resulting from a Dirac function input, $w(t)$.  Moreover, in the multi-input/multi-output (MIMO) case, the $\cH_2$-norm provides a uniform bound on the time response magnitude assuming a disturbance with unit $L_2$-energy:
$\displaystyle \max_{t\geq 0} \|z(t)\|_{\infty} \leq \left\| F( \cdot \,;g)  \right\|_{\cH_2}$ when $\|w\|_{L_2}\leq 1$.

In order to minimize uniformly the output $z$ under the influence of the disturbance, $w$,  we will proceed by determining the ``best" damping such that $\left\| F  \right\|_{\cH_2}$ is minimal.   Using standard theory (see, e.g., \cite{Blanchini12, Burl1998, Zhou96}), it can be shown that the $\cH_2$-norm of the system, $\left\| F  \right\|_{\cH_2}$,  can be expressed via the solution of a Lyapunov equation:
 \begin{equation}\label{transferFunction}
\left\| F (\cdot\,;g) \right\|_{\cH_2}=\left (\frac{1}{2\pi} \trace{\widehat{E}^T \widehat{\mathbb{X}} \widehat{E}}\right )^{\frac{1}{2}} ,
\end{equation}
 where $\widehat{\mathbb{X}}$ solves
\begin{equation}
 \quad \widehat{A}^T(g) \widehat{\mathbb{X}} + \widehat{\mathbb{X}} \widehat{A}(g) =-\widehat{H}^T \widehat{H} \label{Lyap eq},
\end{equation}
and the matrices $\widehat{A}(g)$, $\widehat{E}$, and $\widehat{H}$ are given in (\ref{matrixAHE}).

Although this makes evaluation of $\left\| F(\cdot\,;g)  \right\|_{\cH_2}$ amenable to numerical computation,  the computational resources required to approach realistic problems is still substantial.  Moreover, one observes that the objective function, $\left\| F(\cdot\,;g)  \right\|_{\cH_2}$ will have many local minima with respect to damping positions as encoded in $B$, and for each $B$ there may be many local minimizers with respect to the damping gains, $g=(g_1, g_2, \ldots, g_p)$.  Not surprizingly, many function evaluations are necessary to carry out optimization
with respect to both $g$ and $B$ and this frequently creates an unmanageable computational burden.

 We introduce an approach here which calculates reduced second-order systems in such a way that allows efficient approximation of the
 $\cH_2$ norm,  which in turn brings us a significant  acceleration of the optimization process.  We focus first on efficient optimization with respect to $g$, and then optimization with respect to damping positions; both can be well approximated and cheaply obtained using a reduced order model.  Since we are dealing with structured second-order systems, we use structure preserving methods which are derived particularly for a parametric setting.

There are several  different methods for calculating a reduced system for second-order settings; see, e.g.,\cite{bonin2016fully,Fehr2013,bai2005drs,meyer1996balancing,su1991model,reis2008balanced,breiten2016structure,coprime,BeatGuger05}.  A review of different methods of dimension
reduction, both parametric and nonparametric, can be found in \cite{BennerGW15, ANT05, BAIDEMMEL, BenMS05, Antoulas01asurvey,BauBF14,BenCOW17,hesthaven2016certified}.
The approximation of optimal damping using dimension
reduction for stationary second-order systems was studied in \cite{BennerTomljTruh10, BennerTomljTruh11} where the authors considered optimization of passive damping.   Another approach based on dominant poles, presented first in \cite{SaadvandiMeerbergenDesmet13}, was studied in \cite{BennerKuerTomljTruh15}.
The approach that we present here uses interpolatory projections to produce high fidelity reduced-dimension second-order systems that are then
optimized with respect to damping as measured with the $\cH_2$ system norm.  We employ a variant of the Iterative Rational Krylov Algorithm (IRKA), which is a popular approach for producing high quality reduced models. IRKA produces locally  $\cH_2$-optimal reduced models, dovetailing perfectly with the optimization task at hand.  In Section \ref{IRKAandDampOpt}, we describe the original IRKA approach as well as our variant, and organize its use in damping optimization. Implementation issues are discussed in Section \ref{ImplDetails}. In Section \ref{NumExp}, we describe a variety of numerical experiments that show the advantages of our approach.

\section{IRKA and Damping Optimization} \label{IRKAandDampOpt}

The use of reduced models in parameter optimization involves iterating a two step process wherein one performs parameter optimization on a high-fidelity, parameterized surrogate model in the first, comparatively inexpensive step, which is expected to bring the parameters into the vicinity of their optimal values, followed by a (generally more expensive) update step that corrects deviations that the surrogate model may have with the true parameterized system in the vicinity of the current parameter values.
 This approach is aligned with classical trust region strategies, see for example, \cite{alexandrov1998trf, Gratton2014}.
Since we are seeking to minimize the $\cH_2$ system norm with respect to a variety of damping configurations, an effective choice of surrogate models will be $\cH_2$ optimal reduced order models.
 Efficient interpolatory projection methods have been developed to derive locally $\cH_2$-optimal reduced order models  in related settings (e.g. see \cite{Gugercin_h2model, morBauBBG11, BenBre11,flagg2015multipoint}).

Projection-based methods for model reduction can be described quite simply.  We approximate the full state, $q(t)$, using $q(t)\approx X_r q_r(t)$ where $X_r\in \mathbb{R}^{n\times r}$ has linearly independent columns spanning a \emph{right modeling space} that is still to be determined.  A complementary \emph{left modeling space} spanned by the columns of a second matrix, $W_r\in \mathbb{R}^{n\times r}$, allows us to enforce Petrov-Galerkin conditions that define reduced model dynamics:
\begin{equation} \label{twoSidedProj}
  \begin{array}{c}
  M_r\ddot q_r(t)+C_r\dot q_r(t)+K_r q_r(t)=E_r w(t),\quad \mbox{where}\\[3mm]
   M_r=W_r^* MX_r, \   C_r(g)=W_r^*C(g)X_r, \
K_r=W_r^*KX_r ,\ \mbox{and} \ E_r =W_r^*E.
\end{array}
\end{equation}
The reduced model output then appears as $z_r(t) =H_r q_r(t)$ with $H_r =H X_r$.  Evidently, one should choose $X_r$ and $W_r$, or equivalently, right and left modeling subspaces, so as to ensure $z_r(t)\approx z(t)$ over a wide range of inputs, $w(t)$.
 In \cite{BennerKuerTomljTruh15}, the authors have chosen the columns of $X_r$ and $W_r$ to contain eigenvectors of the polynomial pencil, $\lambda^2 M+\lambda C(g) +K$, that correspond to dominant poles of the transfer function, $F(s; \,g)$, defined in (\ref{TF}).   With this approach one maintains the dominant terms from $F(s; \,g)$ in the reduced transfer function,
\begin{equation} \label{red2ndOrdModel}
F_{2r}(s;\, g)=H_r (s^2M_r+sC_r(g)+K_r)^{-1} E_r.
\end{equation}
Note that $F_{2r}(s;\, g)$  is an $m_{out}\times m_{in}$ complex matrix (same size as $F(s;\, g)$) but now with rational functions of lower order $(2r-1,2r)$ as elements.

We will make a different choice for $W_r$ and $X_r$, choosing them instead so as to enforce \emph{tangential interpolation} conditions:   for selected interpolation points $\sigma_1,\sigma_2,\ldots,\sigma_r\in \mathbb{C}$ and directions $b_1,\ldots, b_r$ and $c_1,\ldots, c_r$,
 we will choose $W_r$ and $X_r$ so that
$$  c_i^T F(\sigma_i) = c_i^T F_{2r}(\sigma_i),\quad
  F(\sigma_i)b_i  =  F_{2r}(\sigma_i) b_i,\quad \mbox{and}\quad
  c_i^T F'(\sigma_i)b_i  =  c_i^T F'_r(\sigma_i) b_i,
$$
for $i=1,\ldots, r$.
Calculation of $W_r$ and $X_r$ that enforce these  interpolation properties is straightforward, and depending on the context, a variety of choices for interpolation points and tangent directions could be used, see, e.g.,
\cite{gallivan2005model,BenBreDam11,Gugercin_h2model,Ant2010imr,BeaG17, BunseGerstner2010}.
This is the thrust of \emph{interpolatory projection methods} for model reduction.  Since we are dealing with structured second-order mechanical systems, we will constrain the choice of $W_r$ and $X_r$ so as to preserve second-order structure as well as important system properties such as system stability, passivity, and reciprocity which is encoded in system structure through the symmetry and positive-definiteness of $M$,  $C$, and  $K$.   This can be accomplished by setting $W_r=X_r$, and then (\ref{twoSidedProj}) becomes:
\begin{equation} \label{oneSidedProj}
 \begin{array}{c}
  M_r\ddot q_r(t)+C_r\dot q_r(t)+K_r q_r(t)=E_r w(t),\quad \mbox{where}\\[3mm]
 M_r=X_r^* MX_r, \   C_r(g)=X_r^*C(g)X_r, \ K_r=X_r^*KX_r, \ \mbox{and} \ E_r =X_r^*E.
 \end{array}
\end{equation}

Using either (\ref{twoSidedProj}) or (\ref{oneSidedProj}), one may express the reduced
transfer function $F_{2r}(s;\, g)$ in (\ref{red2ndOrdModel}) in terms of its $2r$ poles and residues:
\begin{equation} \label{redsysPoleResForm}
F_{2r}(s) =\sum_{k=1}^{2r}\ \frac{c_k\,b_k^T}{s-\lambda_k}
\end{equation}

 First-order necessary conditions for an order $2r$ rational function, $\widehat{F}_{2r}(s)$, having the
 form (\ref{redsysPoleResForm}) (and having $2r$ distinct poles, $\{\lambda_k\}$)  to be an \emph{optimal} $\cH_2$ reduced order approximation to $F(s)$ are due to Meier
 and Luenberger \cite{meier1967approximation}.  They require that $\widehat{F}_{2r}(s)$ be a Hermite interpolant to the full-order system at points in the complex plane that reflect the reduced system poles
 across the imaginary axis and for MIMO systems this need only happen in particular directions in the input/output spaces see, e.g., \cite{Gugercin_h2model,vandooren2008hom,BunseGerstner2010}.  Specifically, if
$\displaystyle \widehat{F}_{2r}(s) =\sum_{k=1}^{2r}\ \frac{\hat{c}_k\, \hat{b}_k^T}{s-\hat{\lambda}_k}$
and $\widehat{F}_{2r}$ is a local minimizer of $\left\| \widehat{F}_{2r} - F \right\|_{\cH_2}$,  then
$\widehat{F}_{2r}$ is also a tangential Hermite interpolant of  $F(s)$
 at $-\hat{\lambda}_k$, $k=1,\ldots,2r$ in the sense that
\begin{equation}  \label{firstOrderOptCond}
\begin{array}{c}
\widehat{F}_{2r}(-\hat{\lambda}_k)\hat{b}_k=F(-\hat{\lambda}_k)\hat{b}_k,\quad
\hat{c}_k^T\widehat{F}_{2r}(-\hat{\lambda}_k)=\hat{c}_k^TF(-\hat{\lambda}_k), \\[3mm]
\quad\mbox{and} \quad
 \hat{c}_k^T\widehat{F}_{2r}'(-\hat{\lambda}_k)\hat{b}_k=\hat{c}_k^TF'(-\hat{\lambda}_k)\hat{b}_k,
\quad \mbox{for }k=1,\ldots,2r.
\end{array}
\end{equation}
IRKA \cite{Gugercin_h2model} is an algorithm that can produce such an $\widehat{F}_{2r}$ efficiently and will provide directly a standard realization for it:
\begin{equation}  \label{redsystemRealization}
\begin{array}{l}
\dot x_{2r}(t)=\widehat{A}_{2r} x_{2r}(t)+\widehat{E}_{2r}w(t),   \\
z_{2r}(t)=\widehat{H}_{2r}x_{2r}(t)
\end{array}
\end{equation}
so that $\widehat{F}_{2r}=\widehat{H}_{2r}\left( s\, I - \widehat{A}_{2r} \right)^{-1}\widehat{E}_{2r}$, but significantly for our application,
$\widehat{F}_{2r}$ \emph{will not} generally have the form of a second-order system transfer function such as $F_{2r}$ in (\ref{red2ndOrdModel}), and this is what motivates the modification we propose for IRKA.

Assume that we have a damping configuration given by $B$ and $G=\diag{(g_1 , g_2 , \ldots, g_p)}$, and we
 would like to obtain a reduced model, $F_{2r}(s;g)$, of the form (\ref{red2ndOrdModel}) that will accurately represent $F(s;g)$ at least for small changes that may be made to the gains, $g_1 , g_2 , \ldots, g_p$.     Using $\|F_{2r}(s;g)\|_{\cH_2}$ as a surrogate for $\|F(s;g)\|_{\cH_2}$, we choose $g_1 , g_2 , \ldots, g_p$ so as to minimize the value of $\|F_{2r}(s;g)\|_{\cH_2}$ (and hopefully, by proxy, $\|F(s;g)\|_{\cH_2}$).

Algorithm \ref{IRKAalg}, which we will also refer to as \textsf{sym2IRKA}, calculates a modeling basis $X_r\in \mathbb{C}^{n\times r}$
using a one-sided projection approach inspired by IRKA.  In each iteration, we form a $2r$ order reduced system transfer function given by \eqref{red2ndOrdModel} and \eqref{oneSidedProj} that interpolates the true transfer function at a set of $r$ interpolation points in $r$ directions
that had been determined in the previous step.   In order to proceed, in Step  \ref{step reduction IRKAalg}, the $2r$ order reduced system transfer function is further reduced to an order $r$ transfer function represented in pole-residue form as
$\widetilde{F}_r(s) =\sum_{k=1}^{r}\ \frac{\tilde{c}_k\,\tilde{b}_k^T}{s-\mu_k}$.  By reflection across the imaginary axis, this gives $r$ interpolation points and tangent directions with which to continue to the next step.  There are many strategies with which to carry out this ``internal reduction" step,  and in what follows we consider three different strategies:
\begin{itemize}
  \item [(a)] \textbf{internal reduction based on balanced truncation}: we use the \emph{balanced truncation method} applied to a linearized order $2r$ realization of $F_{2r}$, to produce an order $r$ (standard system) realization.  The $r$ poles of this realization are reflected to produce the next set of interpolation points.  Details regarding the balanced truncation method can be found, e.g., in \cite{BennerGW15, ANT05, Antoulas01asurvey}.
  \item [(b)] \textbf{internal reduction based on IRKA}: We follow a process similar to (a) except we use a variant of IRKA to produce an order $r$ (standard system) realization.  A natural way of doing this would be to use the original formulation of IRKA, which, upon convergence, will produce a reduced system satisfying necessary conditions for $\mathcal{H}_2$ optimality.  Interestingly, we found this approach did not perform as well as  a ($\mathcal{H}_2$ suboptimal) one-sided modification of IRKA similar to Algorithm 1 (i.e., using $W_r=X_r$) but applied instead to a linearized order $2r$ realization of $F_{2r}$ in order to reduce it to a standard order $r$ realization.    This modified approach appears to reduce the potential for a loss of stability at non-interpolating damping configurations and we generally observe faster convergence of the iteration.
 \item [(c)] \textbf{internal reduction based on dominant poles}:  we choose the $r$ most dominant poles that are closed under conjugation, see e.g., \cite{SaadvandiMeerbergenDesmet13, BennerKuerTomljTruh15});
the transfer function \eqref{TF} can be represented as
$\displaystyle F(s)=\sum_{i=1}^{2n} \frac{R_i}{s-\lambda_i}$ with residues
$
R_i=(Hx_i)(y_i^*E)\lambda_i\in\C^{s\times m},
$
 where   $\lambda_i\in \mathbb{C}$ and $x_i,~y_i\in\Cn\backslash\lbrace0\rbrace$ are,
respectively, eigenvalues, right eigenvectors, and left eigenvectors of the quadratic
eigenvalue problem
\begin{align}\label{qep}
(\lambda_i^2 M+\lambda_i C+K)x_i=0,\quad
y_i^*(\lambda_i^2 M+\lambda_i C+K)=0,\quad i=1,\ldots,2n.
\end{align}
 Although there are a variety of definitions for ``dominant" poles; we take as dominant those poles producing the largest values of
$\frac{\|R_i\|}{|\re{(\lambda_i})|}$.
This choice has been shown to have good performance within the optimization setting that we consider here.
If the $r$ dominant poles fail to be closed under conjugation then we consider taking $r+1$ dominant poles, so as to be closed under conjugation.
\end{itemize}

\begin{algorithm}[h!] 
 \caption{\!:\,\textsf{sym2IRKA}. Symmetrized IRKA for second-order systems}   
\label{IRKAalg}                                            
\begin{algorithmic} [1]                                        
 \REQUIRE  System matrices defining \eqref{MDK} with (fixed) gains $\hat{g}=(g_1 , g_2 , \ldots, g_p)^T$.\\ 
  $\mathsf{itMax}$ - maximum number of iterations allowed for system reduction\\
   Initial shift selection:  $\{\sigma_1\ldots,\sigma_r\}$  (closed under
conjugation). \\
Initial tangent directions: $\{b_1,\ldots, b_r\} $  (also closed under conjugation).
\ENSURE  Modeling basis $X_r(\hat{g})$ producing an interpolatory ROM at $g= \hat{g}$.
   \FOR{$j=1,\ldots, \mathsf{itMax}$}
   \STATE $X_r = [ (\sigma_1^2 M+\sigma_1C+K)^{-1} E b_1,\ldots, (\sigma_r^2 M+\sigma_rC+K)^{-1} E b_r];$
   \STATE  Form reduced system determined by    \\
   $M_r=X_r^* MX_r,C_r=X_r^*CX_r$, $K_r=X_r^*KX_r$,
 $E_r=X_r^*E$ and $H_r=H X_r$
   \STATE Internal reduction step using strategies (a), (b), or (c): \\
   Reduce the order $2r$ transfer function, $F_{2r}$, defined in \eqref{redsysPoleResForm}
   to an order $r$ transfer function, $\widetilde{F}_r(s) =\sum_{k=1}^{r}\ \frac{\tilde{c}_k\,\tilde{b}_k^T}{s-\mu_k}$,
 such that $\mu_1,\ldots, \mu_r$,  are closed under conjugation \label{step reduction IRKAalg}
   \STATE Assign $\sigma_k = -\mu_k$ and $b_k= \tilde{b}_k$ for $k=1,\ldots,r$
      \IF  {$\{\sigma_k\}$ converged}\label{conver_test}
   \STATE    \textbf{break};
   \ENDIF
   \ENDFOR
\STATE $X_r(\hat{g})\leftarrow X_r$;
\end{algorithmic}
\end{algorithm}
Other potential strategies for internal reduction may be found in \cite{Wyatt12}.  Note that  the  \textbf{for} loop in \textsf{sym2IRKA}  goes up to a specified maximum number of iterations, $\mathsf{itMax}$, but is stopped earlier if two consecutive sets of interpolation points, $\{\sigma_i\}$ have changed only slightly within a given tolerance. In terms of computational effort, note that Step 4 deals only with reduced systems, requiring effort that scales with $r\ll n$ rather than $n$.

Although we are solving large linear systems repeatedly within Step 1 and Step 6, we find that these tasks can be implemented efficiently through the  use of the modal coordinates. This is discussed in the next subsection.

Through the use of \textsf{sym2IRKA}, we obtain a reduced model that has the same structure as the true system and should replicate the response of the true system (and hence replicate its $\mathcal{H}_2$ norm) to relatively high accuracy at the given sampling gain  $(g_1 , g_2 , \ldots, g_p)$.  In order to approximate  the $\mathcal{H}_2$ norm of $F(s,g)$ with the (easily computed)
$\mathcal{H}_2$ norm of $F_{2r}(s,g)$ over a sufficiently broad range of parameter values in $g$, we will employ a parametric model order reduction approach described in \cite{morBauBBG11}. In this approach we calculate
$\mathcal{H}_2$-based interpolatory reduced systems for a small number of selected sampling parameters and then form an aggregate interpolatory reduced model. The original and reduced model then have virtually the same $\mathcal{H}_2$ norms on the given set of sampling gains and generally nearly so at points in between. This makes the optimization process more robust and damping optimization can be performed more efficiently.  A similar approach was also used in \cite{BennerKuerTomljTruh15} where the authors considered local approximations based on dominant poles of the system.

In our setting, parameterized model order reduction (PMOR) is organized so that we calculate reduced systems given by \eqref{oneSidedProj} for several, say $m$, damping gain configurations, $g^{(i)}$, $i=1,\ldots, m$ and for each of these sampled gain configurations (represented as a $p$-tuple of individual damper gains),  we calculate a modeling basis, $X^{(i)}_r$, that corresponds to a reduced system given by \eqref{oneSidedProj}.  We then merge all $X^{(i)}_r$ $i=1,\ldots, m$ into an aggregate modeling basis, $X=[\, X^{(1)}_r,\, X^{(2)}_r,\,\ldots,\,X^{(m)}_r\,]$ and form an aggregate reduced system model,
 \begin{align}\label{gain_redsys}
\begin{split}
 X^* MX\ddot q_k(t)+X^*C(g)X\dot q_k(t)+X^*KXq_k(t)&=X^*E w(t),\\
z(t)&=H Xq_k(t).
\end{split}
\end{align}
  By using this approach, we are assured of a good approximation of system response (and $\mathcal{H}_2$ norm approximation) for any of the particular gain configurations, $g^{(i)}$, $i=1,\ldots, m$, and we will generally retain good approximation for nearby gain configurations as well.  Since we are interested in optimizing these gains, we may organize the procedure so that sampled gain configurations are determined adaptively during the optimization process, and used to augment the aggregate modeling basis.  Gain configurations that were useful earlier in the process but that no longer contribute significant information can be dropped from the aggregate modeling basis.   How best to produce a balanced parameter sampling strategy is an important and generally unresolved issue, but for particular problems there have been some advances (e.g., see the recent survey \cite{BennerGW15}).
For our damping optimization problem, we consider two main strategies to sample damping gain configurations:
\begin{itemize}
  \item[(i)] \textbf{predetermined sampling of damping gain configurations}:
  We use a predetermined set of $m$ damping gain configurations: $g^{(1)},\ldots,g^{(m)}$
that has been chosen from the set of feasible damping gain configurations, $g^{(k)}\in \bigtimes_{i=1}^p [g_i^-,g_i^+]$.
This choice can be done via uniform sampling across a fixed mesh in the feasible region,
or it can include \emph{ad hoc} choices of damping gain configurations as well.
For example, we include among our predetermined damping gain configurations,
the trivial configuration, $g=0$, that is, we include a system having only internal
damping and no external damping.  This approach
is described in more detail as Algorithm \ref{OptGain-i}.

 \item [(ii)] \textbf{adaptive sampling during optimization}.  Starting with an initial damping gain configuration (say, $g^{(1)}=0$),
  construct a reduced order model, $F_{2r}(s;\, g^{(1)})$ using \textsf{sym2IRKA}.  Then, allowing $g$ to vary, determine the next
  damping gain configuration to be sampled by finding $g^{(2)}$ that solves (approximately)
  \begin{equation*}
g^{(2)}=\argmin_g \left\| F_{2r} (\cdot\,;g) \right\|_{\cH_2}.
\end{equation*}
Repeat this process, each time augmenting the modeling subspace used to construct reduced models with information from
the newest reduced order model.  This is discussed below and described in more detail as Algorithm \ref{OptGain-ii}.
 \end{itemize}

\begin{algorithm}[ht] 
 \caption{Approximation of optimal gains using predetermined gain configuration samples}   
\label{OptGain-i}                                            
\begin{algorithmic} [1]                                        
\REQUIRE  System matrices defining \eqref{MDK}; A set of $m$ (different) gain configurations $g^{(1)},\ldots,g^{(m)}$.\\
The number of retained poles $\mathsf{nRetPoles}$ for each gain configuration;\\
 Initial choices for shift selection $\sigma_1\ldots,\sigma_r$ and directions $b_1,\ldots, b_r$; \\
 \ENSURE  Approximate optimal gains.                                            
\FOR {$j=1,\ldots,m$}
\STATE  Using the gain configuration $g^{(j)}$, calculate a reduced order modeling subspace, $V^j$, with \textsf{sym2IRKA} (Algorithm \ref{IRKAalg}).
\ENDFOR
\STATE  Aggregate modeling spaces, into $X=\orth( [ V^1,\ldots, V^m ])$.
\STATE Form a  global reduced system using $X$ as in \eqref{gain_redsys}.
\STATE Find (approximate) optimal gains by optimizing the gains of the global reduced system, \eqref{gain_redsys} (using an appropriate optimization tool).
\end{algorithmic}
\end{algorithm}

\begin{algorithm}[!h] 
 \caption{Computation of optimal gains with adaptive sampling}   
\label{OptGain-ii}                                            
\begin{algorithmic} [1]                                        
\REQUIRE  System matrices defining \eqref{MDK}; An initial gain configuration $\hat{g}^{(0)}=(g^{(0)}_1 , g^{(0)}_2 , \ldots, g^{(0)}_p)^T$.\\
 Initial choices for shift selection $\sigma_1\ldots,\sigma_r$ and directions $b_1,\ldots, b_r$; \\
 $\mathsf{nRetPoles}$ - number of retained poles for each gain configuration;\\
  $\mathsf{itMax}$ - maximum number of iterations for \textsf{sym2IRKA} (Algorithm \ref{IRKAalg}).\\
    $\mathsf{tolDiff}$ - termination criterion for gain optimization  \\
 \ENSURE  Approximate optimal gains.
 \STATE $j=0$;
\REPEAT
\STATE  Using the gain configuration $g^{(j)}$, calculate a reduced order modeling subspace,
                $V^j$, with \textsf{sym2IRKA}. \label{IRKA call in ii}
\STATE Form a reduced system using  $X=\orth( [V^0,V^1,\ldots, V^j ])$ as in \eqref{gain_redsys}.
\STATE $j=j+1$
\STATE Find an (approximately) optimal gain configuration by optimizing the gains of the global
reduced system, \eqref{gain_redsys} (using an appropriate optimization tool), and denote it by $g^{(j)}$
\UNTIL {$|g^{(j)}-g^{(j-1)}|<\mathsf{tolDiff}$}
\STATE return $g^{(j)}$
\end{algorithmic}
\end{algorithm}

Strategies (i) and (ii) (and the associated Algorithms \ref{OptGain-i} and \ref{OptGain-ii}) could be viewed as representing two extremes and one can easily construct sensible sampling strategies that combine elements of both in various ways.   The common element in both strategies is the use of a reduced order model that inherits from the original model, \eqref{TF}, the parameterization with respect to damping.  The optimal damping configuration then is sought by optimizing the reduced order model, using it as an inexpensive surrogate for the original model.
This is effective when the $\cH_2$-norm of the reduced order model closely tracks  the $\cH_2$-norm of the original model as the damping parameters in $g$ change.   In order to assure that, one may include in the optimization goals a \emph{greedy search} for the damping
 parameter $\hat g$, such that the deviation of $\cH_2$-norm values between full-order and reduced-order models is largest:
\begin{equation}\label{errorTF}
\hat g=\argmax_g |\left\| F (\cdot\,;g) \right\|_{\cH_2}-\left\| F_{2r} (\cdot\,;g) \right\|_{\cH_2}|,
\end{equation}
where, as before, $F$ is transfer function of full system while $F_{2r}$ is transfer function of reduced system given by \eqref{gain_redsys}.
Augmenting the current reduced model so as to interpolate the full order model at  $\hat g$ is likely to improve fidelity over a wide range of intervening $g$.   One may repeat the procedure  until the deviation given by \eqref{errorTF} is acceptable.  Such an approach ensures that   the reduced system is an accurate surrogate for the full system viewed as a function of damping parameters, $g$.

There are two principal disadvantages to this approach in our setting. The first is that the determination of $\hat g$ from \eqref{errorTF} requires evaluation of the $\cH_2$-norm of the full system, which either is computationally very demanding, or if the $\cH_2$-norm is computed only approximately, this introduces an additional layer of approximation to the approach.  The second disadvantage is that it may happen that the greedy selection of gains produced by \eqref{errorTF} may be far away from the gains that minimize the $\cH_2$-norm of the full system, so we expend significant effort in producing a reduced model that has high fidelity at damping configurations that are not interesting for us.

Therefore, we propose an adaptive sampling strategy that samples gain configurations that have approximately minimized the $\cH_2$-norm of an intermediate reduced order system.   Assume that we have obtained at some stage in the process a reduced system given by \eqref{gain_redsys} with $X=V^1$.  We then determine a gain configuration, $\hat g$, that minimizes the surrogate objective function:
\begin{equation*}
\hat g=\argmin_g \left\| F_{2r} (\cdot\,;g) \right\|_{\cH_2},
\end{equation*}
where  $F_{2r}$ is the transfer function of the reduced system given by   \eqref{gain_redsys}.   We generate another modeling subspace, $V^2$, that would produce a high fidelity reduced model at the system with a damping configuration of $\hat g$, (i.e., we calculate $V^2$ using Algorithm \ref{IRKAalg} for the gain configuration $\hat g$), but we update the reduced system by augmenting the modeling subspace with $V^2$, instead of replacing it with $V^2$: $X=\orth([V^1,\, V^2])$. These steps are repeated (and the modeling subspace grows), until the difference between consecutively sampled gains drops below the prescribed tolerance $\mathsf{tolDiff}$.
%
%
until the difference between consecutive sampling gains is equal or smaller than prescribed tolerance $\mathsf{tolDiff}$.
In this way, the modeling subspace spanned by the columns of $X$ continues to grow and the reduced order model generally provides progressively higher fidelity approximations to the full order system in the vicinity of the sampled damping configurations.  However, this added fidelity might not provide improved information about the optimum damping configuration.   Note also that if too low an order is used for the reduced model then it may have poor fidelity notwithstanding its local $\cH_2$-optimality and the reduced models could fail to recover even coarse $\cH_2$-norm information for the full system resulting (potentially) in a false optimum.  As a practical matter, we have found that quite small reduced system orders will still produce reasonable $\cH_2$-norm estimates and we have not observed such failures.
It is significant that we are able to avoid entirely the evaluation of the
$\cH_2$ norm of the full system, which is a forbidding computational challenge for large scale systems.

%

\section{Some Implementation Details}  \label{ImplDetails}
%
%

Our approach relies on the effectiveness of \textsf{sym2IRKA} (Algorithm \ref{IRKAalg}), in generating structured reduced order models with good fidelity.  The efficiency of this process is improved if we introduce \emph{modal coordinates}:  Since $M$ and $K$ are symmetric positive definite matrices there exists a nonsingular matrix $\Phi$ (the \emph{modal matrix}) which simultaneously diagonalizes both $M$ and $K$:
\begin{equation}\label{simult. diag of M,K}
\Phi^T K \Phi = \Omega^2= \diag(\omega_1^2, \ldots, \omega_n^2)
\quad \text{and} \quad \Phi^T M \Phi = I,
\end{equation}
where $0<\omega_1\leq\omega_2\leq\ldots\leq\omega_n$ are the undamped natural frequencies of the system  while
the columns of the matrix $\Phi$ are eigenvectors (\emph{modes}) of the undamped system.
We have adopted the usual assumption on internal damping, namely that internal damping is \emph{modal damping}, meaning that $\Phi$ diagonalizes the internal damping matrix, $C_{int}$, as well.  We have taken internal damping to be modeled as a small multiple of critical damping, so $\Phi^TC_{int}\Phi=2\alpha_c \Omega$ with small $\alpha_c$. Evidently, this approach can be adapted immediately to other  models of internal damping as well, e.g.,  Rayleigh damping or more general proportional damping (see, e.g. \cite{KuzmTomljTruh12}).

Now, rewrite the system \eqref{MDK} in modal coordinates.
By using \eqref{simult. diag of M,K} and substituting $q(t)=\Phi \hat q(t)$, we
obtain:
\begin{align}\label{MDK1ps}
 \ddot {\hat q}(t)+(\alpha\Omega+\Phi^TBGB^T\Phi)\dot {\hat
q}(t)+\Omega^2{\hat q}(t)&=\Phi^T E w(t),\\
z(t)&=H\Phi \hat q(t).\label{MDK3ps}
\end{align}

We are focused in optimizing parameters in matrix $G=\diag{(g_1, g_2, \ldots, g_p)}$ where number of dampers is usually much smaller than the number of states. This means that within the damping matrix during the optimization process we have small rank update which one should use in order to use the structure of system matrices efficiently. In particular, in steps 1 and 6 of Algorithm \ref{IRKAalg}, one drawback from the computational point of view, is that we need to solve many linear systems especially within the main loop of the algorithm. Here, since we are in modal coordinates we use the following approach.

In applying \textsf{sym2IRKA}, it is necessary to solve repeatedly systems of linear equations having the form
$(\sigma^2 M+\sigma C+K)^{-1} B b$, for varying shifts $\sigma$ and directions $b$.
In modal coordinates,
{\small
\begin{align*}
(\sigma^2 M+\sigma C+K)^{-1} B b & =
  \Phi(\sigma^2 I+\sigma \alpha_c\Omega+\sigma\Phi^TBGB^T\Phi +\Omega^2)^{-1} \Phi^TBb\\
   &= \Phi\mathbb{D}(\sigma)^{-1}\Phi^TBb -\sigma \Phi\mathbb{D}(\sigma)^{-1}\Phi^T B
   \mathbb{M}(\sigma,g) B^T \Phi \mathbb{D}(\sigma)^{-1}\Phi^TBb
\end{align*} }
where
\begin{align*}
  \mathbb{D}(\sigma)  & =\sigma^2 I+\sigma \alpha_c\Omega  +\Omega^2,\quad   \mathbb{G}(g)=\diag(\sqrt{g_1},\ldots,\sqrt{g_p}), \\
 & \mathbb{M}(\sigma,g) =\mathbb{G}(g)(I_p+\sigma \mathbb{G}(g)B^T \mathbb{D}(\sigma)^{-1}B\mathbb{G}(g))^{-1}\mathbb{G}(g),
\end{align*}
  and we have used the Sherman-Morrison-Woodbury formula \cite{GVL89} in passing from the first to the second line.

Note that $\mathbb{D}(\sigma)$ is an $n\times n$ diagonal matrix and
$\mathbb{M}(\sigma,g)$ is a $p\times p$ matrix.  Recalling that $p\ll n$ is the number of dampers,
the evaluation of products of the form $\mathbb{M}(\sigma,g)\,z$ or even outright evaluation of the full matrix $\mathbb{M}(\sigma,g)$
can be implemented very efficiently.

 The choice of initial gain in Algorithm \ref{OptGain-ii}  follows \cite{BennerKuerTomljTruh15}, where it was found effective to take
 initial gains uniformly equal to zero ($g^{(0)}=0$), so that the optimization process starts with a modally damped system having each mode
 with a presumed fixed initial fraction, $\alpha_c$, of critical damping. In \cite{BennerKuerTomljTruh15}, modal truncation was used to form reduced order models to be used as surrogates in the damping optimization process (see also, \cite{BennerTomljTruh11,Gaw, VES2011}); here we improve on the standard modal approximation by using IRKA approach on system without external damping.
It is shown that one can take undamped eigenvectors that corresponds to the smallest $r$ undamped eigenfrequencies, in order to obtain modal approximation of the system. This improves robustness of optimization process and it can also be calculated off-line since it does not include information on external damping. Here we will improve standard modal approximation (see, e.g. \cite{BennerTomljTruh11, BennerKuerTomljTruh15, Gaw, VES2011}) by using \textsf{sym2IRKA} (Algorithm \ref{IRKAalg}) without external damping (g=0).

The application of Krylov-based model reduction of second-order systems with proportional damping was studied in \cite{BeatGuger05}, but since in this particular case (with $g=0$) in Step 1 and Step 6 of \textsf{sym2IRKA}, we need only solve diagonal systems and running \textsf{sym2IRKA} can be done off-line without significant cost. This first approximation determines an initial shift selection $\sigma_1\ldots,\sigma_r$.   The corresponding  tangent directions $b_1,\ldots, b_r$ are calculated as right singular vectors of reduced system transfer function evaluations on the corresponding initial shifts.  In this way, the first cycle of either Algorithm \ref{OptGain-i} or \ref{OptGain-ii} involving the development of a reduced order model for the case $g=0$ can be counted as part of the off-line phase.
We have observed in the numerical examples described below, that this initial model reduction approach provides ultimately a better approximation of optimal gains  than if we were to use modal approximation for the initial surrogate model.

%

Note that  in Step \ref{IRKA call in ii} of  Algorithm \ref{OptGain-ii}, we recycle information obtained from earlier \textsf{sym2IRKA} runs:  optimal shifts and directions obtained from \textsf{sym2IRKA} can be used as starting shifts and directions for the next run of \textsf{sym2IRKA}. This  further  improves the convergence behaviour of \textsf{sym2IRKA}, especially when consecutive gain configurations are close to one another.

\section{Numerical experiments} \label{NumExp}

We  illustrate here the advantages that accrue with the use of \textsf{sym2IRKA}, which is our variant of IRKA introduced as Algorithm \ref{IRKAalg}.  We present results obtained using different sampling strategies for intermediate damping configurations and different approaches taken for the internal reduction step (Step \ref{step reduction IRKAalg}) in \textsf{sym2IRKA}.  These results will be compared with the approach presented in \cite{BennerKuerTomljTruh15},  which describes another approach to damping optimization.

\begin{example}\label{ex1}
{\em
We consider a mechanical system of $n$ masses connected in sequence by simple springs
(see \cite{TrTomPuv16}), so that each mass  is connected via a spring to two neighboring masses.
The mathematical model is given as
(\ref{MDK}), where the mass and stiffness matrices are
\begin{align*}
M&=\diag{(m_1,m_2,\ldots,m_n)}  ,\\
K&=\left(%
\begin{array}{ccccccc}
  2k_1+2k_2       & -k_2      &-k_3  &   & && \\
  -k_2            & 2k_2+2k_3 & -k_3 &-k_4  &  &&\\
  -k_3            & -k_3      & 2k_3+2k_4 & -k_4 &-k_5  & & \\
                  & \ddots   & \ddots & \ddots & \ddots &&\\
                  &     &  \ddots &2k_{n-2}+2k_{n-1}& -k_{n-2}    & -k_{n-1} &\\
                  &    &  & -k_{n-1} & 2k_{n-1}+2k_{n} & -k_n  & \\
                  &    &  & -k_{n-2} & -k_n & 2k_n+2k_{n+1} &\\
\end{array}%
\right).
\end{align*}

We consider this system with $n=$1900 masses and the following values for mass and stiffness:
\begin{align*}
    k_i&=500, \quad\mbox{for }i=1,...1900;
  & m_i= \left\{ \begin{array}{ll} 144-\frac{3}{20}i, \quad & i=1,\ldots,475, \\
   \frac{i}{10}+25 , \quad & i=476,\ldots,1900.
                 \end{array}\right.
\end{align*}
The parameter $\alpha_c$ associated with internal damping (\ref{C_int}) is set to $ 0.005$ (i.e., $0.5\%$ critical damping).

The primary excitation comes from 10 disturbances applied at 10 consecutive masses starting with mass 471 until mass 480, in such a way that the magnitude of the excitation is highest at the centerpoint of the range of masses, dropping off at the edges. Thus, the matrix $E\in \mathbb{R}^{1900\times 10}$ is defined so that
 \begin{equation} \label{E2}
  E(471:480,1:10) =\diag(10,20,30,40,50,50,40,30,20,10);
 \end{equation}
all other entries are zero.

For the output, we are interested in observing the displacement history of 18 masses that are uniformly distributed throughout the system starting with mass 100:
 \[
 z(t)=
 \left[q_{100}(t)~ q_{200}(t) ~ \ldots ~q_{1800}(t)\right]^T,
 \]
 which we obtain by choosing matrix   $H\in\mathbb{R}^{18\times 1900}$ as
  \begin{equation*}
 H(1:18,100:100:1800 )=I_{18\times 18};
 \end{equation*}
 all other entries are equal to zero.

External damping is provided by four dampers (i.e., $p=4$) that are configured into two independently positioned pairs of adjacent dampers,
with the damping geometry determined by
\begin{equation}\label{B2}
 B =\left [e_{j} ~  e_{j+1} ~ e_{k}~   e_{k+1} \right],
\end{equation}
where $e_j$ is the $j$th canonical vector and indices $j$ and $k$ determine the damper positions. Since the first and second dampers share the same gain $g_1$, and the third and fourth dampers share the gain $g_2$, $G$ appears as  $$G=\diag{(g_1,~g_1,~g_2,~g_2)}\in \mathbb{R}^{4\times 4}.$$

We consider different damping configurations by varying the indices $k$  and $j$.  In particular, we consider the following
possibilities:
 \begin{equation*}
 \{(j,k)\, |\ j \in \{50,150,250,350\}\ \mbox{and}   \ k\in\{850,950,1050,\ldots, 1850\}\}.\\
 \end{equation*}
 This gives 44 different damping configurations over which we optimize.
 }\end{example}

During the optimization procedure, we use the following parameters in Algorithms \ref{OptGain-i} and \ref{OptGain-ii}:
\begin{align*}
 \quad r=60,\quad \mathsf{tolDiff}=0.001, \quad \mathsf{itMax}=40.
 \end{align*}
A convergence tolerance of $0.001$ was used for termination in step \ref{conver_test} of \textsf{sym2IRKA}.

Viscosities were optimized by   {\sc Matlab}'s built-in \verb"fminsearch" using a termination tolerance of $10^{-4}$,
ending the optimization if either function values or variable values change by less than  and $10^{-4}$.
The starting point for the optimization procedure was $g^{(0)}=(1000,1000)$.



We perform several comparisons with this example. The first comparison is between our approach using \textsf{sym2IRKA} and the dominant pole approach using \cite[Algorithm 1 and 2]{BennerKuerTomljTruh15}.

For clearer comparisons, we use preset gain configurations with  four parameter samples:
$g^{(1)}=(0,0)$ (internal damping only),  $g^{(2)}=(1000,\,1000)$, $g^{(3)}=(100,\,1000)$   and $g^{(4)}=(1000,\,100)$.
When we used adaptive sampling,  then the starting initial gain was taken as $g^{(1)}=(0,0)$ allowing only for internal damping.

We note that the magnitude of the first and second gains varies between 500 and 4000, thus for the sake of better comparisons
with the preset gain configurations, we have chosen values in this range. Usually such information is not known in advance and an adaptive sampling strategy could be advantageous.

\begin{figure}[h!]
\begin{center}\centering
\includegraphics[width=\linewidth]{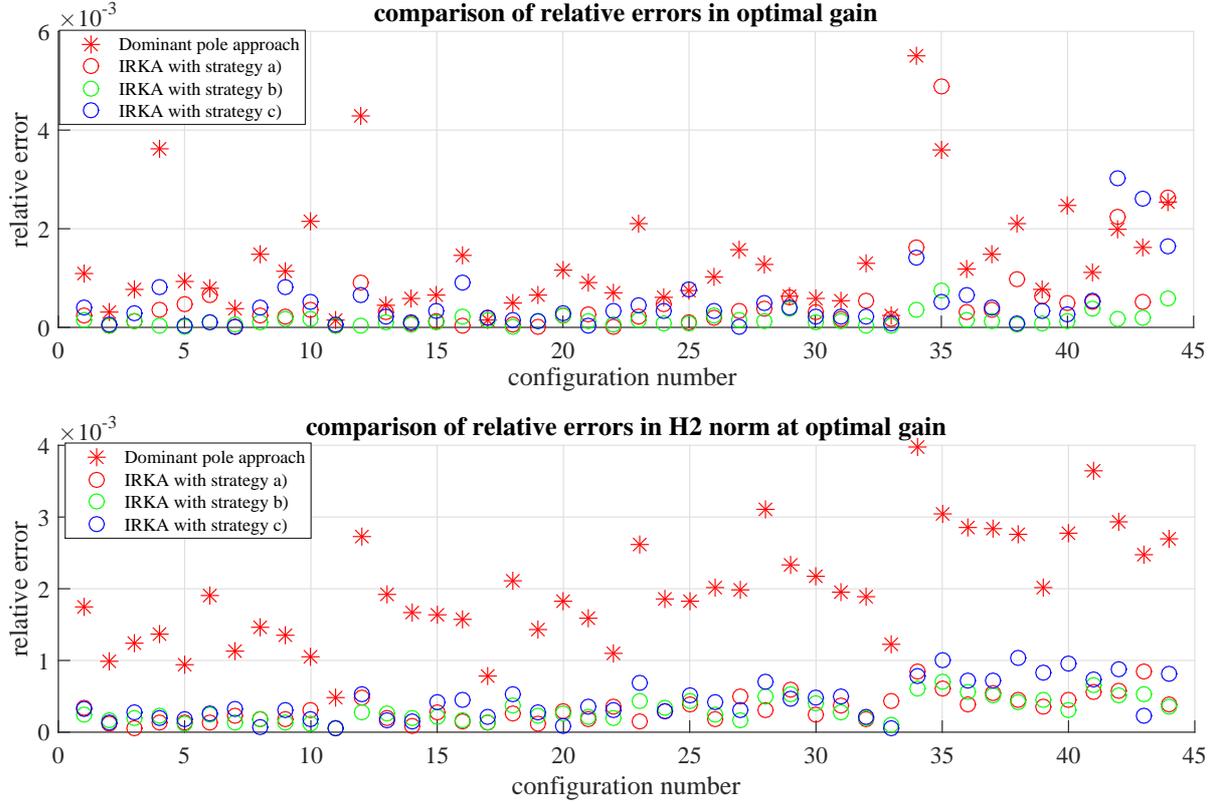}
\caption{Relative errors for the gain and the objective function.}
 \label{rel-error-sDPA}
\end{center}
\end{figure}

As one can see  in the upper plot of   Figure \ref{rel-error-sDPA},  approaches using \textsf{sym2IRKA} have smaller relative errors in optimal gain.   Relative error of the optimal gain was calculated by $\|g^*-g\|/\|g\|$, where $g$ and $g^*$ denote the optimal gain calculated with and without dimension reduction, respectively.  The influence of our $\mathcal{H}_2$-based reduction approach is emphasized even further in the lower plot of  Figure \ref{rel-error-sDPA} where we depict the relative errors in the $\mathcal{H}_2$ cost function.  \textsf{sym2IRKA} consistently yields smaller error.

Next, we consider different damping configuration sampling strategies coupled with different internal reduction strategies ((a), (b) and (c) as described above) used in Step \ref{step reduction IRKAalg} of \textsf{sym2IRKA}.
 \begin{figure}[h!]
\begin{center}\centering
\includegraphics[width=16cm]{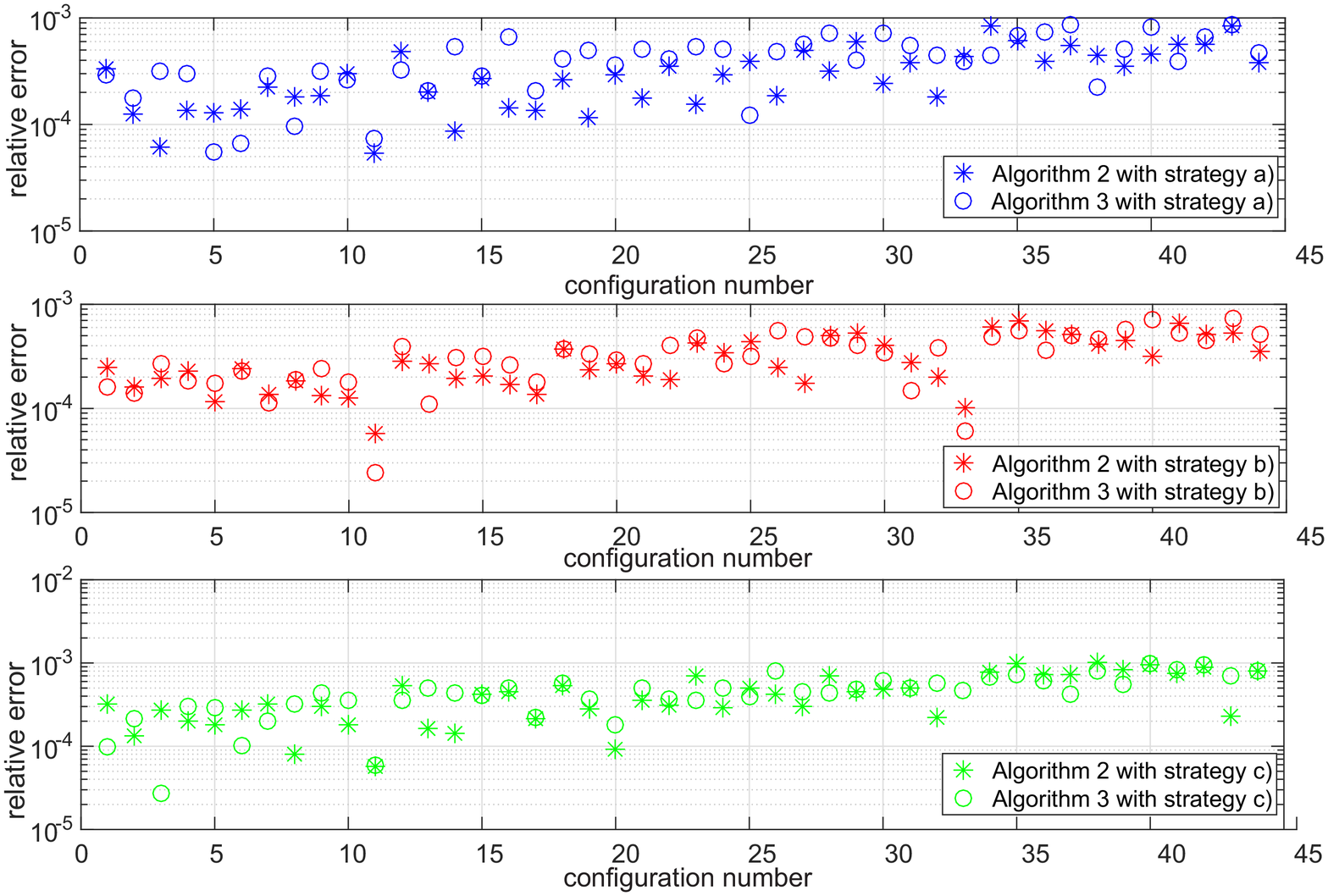}
\caption{Example 1, relative errors for $\cH_2$ norm at optimal gain for Algorithm 2 and Algorithm 3 with strategies a), b) and c)}
 \label{relative-error-abc}
\end{center}
\end{figure}
 For Algorithm \ref{OptGain-i}, a predetermined damping configuration sampling was used. For Algorithm \ref{OptGain-ii} the first initial gain was set to an approximation of optimal gains obtained by using \textsf{sym2IRKA} applied to a system that has only internal damping.
In Figure \ref{relative-error-abc}
we show relative errors for  Algorithm \ref{OptGain-i} and \ref{OptGain-ii} with  internal reduction strategies based on balanced truncation, \textsf{sym2IRKA}, and the dominant pole algorithm related to strategies  (a), (b), and (c) described above.

Besides accuracy, an important comparison criterion is the speed-up of the underlying optimization algorithm. Table 1 shows the average speed-up for the optimization process obtained by \textsf{sym2IRKA} for different strategies of internal reduction. We show here the average ratio between the time  required for the gain optimization with and without the approximation technique using \textsf{sym2IRKA}. Overall, Algorithm 3 yields bigger speed-ups. For example, Algorithm 3 with Strategy (c) has yielded an average speed up of $346.20$.
%
%
 By way of contrast, note that the dominant pole algorithm based on \cite[Algorithm 1 and 2]{BennerKuerTomljTruh15} accelerated the optimization process only by a factor of 43.99 in average cases.

\begin{table}
\begin{center}
\begin{center}
\begin{tabular}{|c||c|c|}
  \hline
     Acceleration factor for:  & Algorithm 2 & Algoritm 3 \\ \hline \hline
   Strategy (a)             &   187.37  & 288.62 \\\hline
   Strategy (b)             &   146.12  & 338.29  \\\hline
   Strategy (c)             &   228.63  & 346.20 \\
  \hline
\end{tabular}
\end{center}
\caption{Time ratios for Example 1,
i.e., acceleration factor of each approach using model reduction accelerated optimization}\end{center}
\end{table}

We conclude that we have obtained satisfactory relative errors in these trials and the approach using \textsf{sym2IRKA} significantly accelerated the optimization process. Moreover, we see that an approach that includes adaptive sampling additionally improves the efficiency of the \textsf{sym2IRKA} approach, especially if feasible intervals of optimal gains are not known in advance.

For the predetermined sampling strategy, we have obtained similar reduced dimension for approaches that use \textsf{sym2IRKA}, however an approach that uses the dominant pole approach had smaller reduced dimension.
In particular, the average reduced dimension for predetermined sampling with strategy (a), strategy (b), and strategy (c) was 150, 178 and 141, respectively, while for the approach based on dominant poles,  we obtained an average reduced dimension of 98.
This means that an approach that uses \textsf{sym2IRKA} provides a modeling subspace that spans a ``richer'' subspace, if we compare with the subspace that is obtained by using the dominant pole approach.  The reduced dimension has a direct impact not only on the time ratio, but also on the rate of convergence of \textsf{sym2IRKA}.

\begin{example}\label{ex2}
{\em
We now consider a different structure that allows for more variables within the optimization process.
Consider a  mass oscillator with $2d+1$ masses and $2d+3$ springs (see  \cite{BennerKuerTomljTruh15, BennerTomljTruh11}). There are two rows of $d$ masses connected with springs; the first row of masses have stiffnesses $k_1$ and the second row have stiffnesses $k_2$. The first masses on the left edge  ($m_1$ and $m_{d+1}$) are connected to  a fixed boundary while  on the other side of rows the masses ($m_d$ and $m_{2d}$)  are connected to   mass $m_{2d+1}$   with a stiffness $k_3$ connected to a fixed boundary.

As in the previous example, the mathematical model is given by (\ref{MDK}), with diagonal mass matrix $M =\diag{(m_1,m_2,\ldots,m_n)}$, but now  the mass and stiffness matrices are given by
\begin{align*}
K =\begin{bmatrix}
    K_{11} &   &  -\kappa_1\\
      & K_{22} &  -\kappa_2\\
     -\kappa_1^T & -\kappa_2^T  &  k_1+k_2+k_3\\
 \end{bmatrix},\quad \mbox{with} \quad
 K_{ii} =k_i
 \begin{bmatrix}
    2 & -1  & & &\\
     -1 & 2 &-1 & & \\
      &     \ddots & \ddots  & \ddots&\\
      &    & -1 & 2 & -1 \\
       & & &-1 & 2\\
 \end{bmatrix},\quad \kappa_i=\begin{bmatrix}0\\
 \vdots\\
 0\\
 k_i
 \end{bmatrix},
\end{align*}
for $i=1,2$.
We will consider the $n$-mass oscillator with $d=1000$, which means that we consider 2001 masses with the following configuration for the masses:
\begin{align*}
    m_{2001}=100,
  & \quad m_i= \left\{ \begin{array}{ll} 100-  \frac{i}{10}, \quad & i=1,\ldots,499, \\
   \frac{i}{30} +33  , \quad & i=501,\ldots,1000, \\
  100  - (i-999)\frac{5}{20}+ \frac{(i-999)^2}{ 5000}  , \quad & i=1001,\ldots,2000,
                 \end{array}\right.
\end{align*}
The  stiffness values are given by
$$ k_1=400, k_2=100 ,k_3=300. $$
The parameter $\alpha_c$ associated with internal damping (\ref{C_int}) is set to $ 0.003$ (i.e., $0.3\%$ critical damping).

The primary excitation are 21 disturbances applied to 21 masses that are closest to ground.  The 10 masses closest to the left-hand side are such that masses closest to the ground are exposed to higher magnitude disturbances.  Additionally, a disturbance is applied to the mass on the right-hand side. That is $E\in \mathbb{R}^{n\times 21}$ with
 \begin{align*}
  E(1:10,1:10) &=\diag(1000:-100:100),\\
    E(1001:1010,11:20) &=\diag(1000:-100:100),\\
    E(2001,21)&=2000;
 \end{align*}
 all other entries are equal to zero.

Regarding the output, we are interested in the 21 displacements of the first row of masses and 21 displacements  of the second row of masses. In particular, we have that
 \[
 z(t)=
 \left[q_{490}(t)~ q_{491}(t) ~ \ldots ~q_{ 510}(t)~ q_{1490}(t)~ q_{1491}(t) ~ \ldots ~q_{ 1510}(t)\right]^T,
 \]
 so that  $H\in\mathbb{R}^{42\times 2001}$.

In this example we illustrate how our approach based on \textsf{sym2IRKA} performs within a setting having more optimization variables. We still
consider four dampers ($p=4$) but now with different viscosities. The geometry of external damping is determined  by a matrix $B$ given by
\begin{equation}\label{B2ex}
 B  =\left [e_{j}-e_{j+5} ~~  e_{j+20}-e_{j+25} ~~ e_{k}-e_{k+5} ~~  e_{k+20}-e_{k+25} \right],
\end{equation}
where $e_j$ is the $j$th canonical vector and indices $j$ and $k$ determine damping positions. Since the dampers now have (potentially) all different viscosities, we have that   $$G=\diag{(g_1,~g_2,~g_3,~g_4)}\in \mathbb{R}^{4\times 4}.$$

The positions are chosen such that the first two dampers are applied on the first row of masses, while the third and the fourth damper  are applied on the second row of masses.  Similarly as in the previous example, in order to consider different damping configurations, we vary the indices $k$  and $j$.
The following
configurations  are considered:
 \begin{equation*}
 \{(j,k)\, |\, j \in \{250,450,650,850\}   \ \mbox{and}   \  k\in\{1150, 1250,  1350, 1450, 1550,  1650,  1750\}\}.\\
 \end{equation*}
 This gives 28 different damping configurations over which we optimize.
 }\end{example}

During optimization, we use the following parameters for Algorithms \ref{OptGain-i} and \ref{OptGain-ii}:
\begin{align*}
 \quad r=120,\quad \mathsf{tolDiff}=0.05, \quad \mathsf{itMax}=40.
 \end{align*}
A tolerance of $0.001$ was used for checking convergence in Step \ref{conver_test} of \textsf{sym2IRKA}.
The starting point for optimization (as performed by \verb|fminsearch|) was set to {\small{$(1000,1000,1000,1000)$}}.
The termination tolerance for \verb|fminsearch| was set to $5\cdot 10^{-4}$.

We observed that the magnitudes of all optimal gains vary from 350 to 7000, thus, as we have done in the previous example, for the predetermined gain sampling we choose values that belong to this range.  In particular, five sampling parameters were used as predetermined samples of damping gain configurations:   $g^{(1)}=\mbox{\small{$(0,0,0,0)$}} $ -- representing a system that has only internal damping, $g^{(2)}=\mbox{\small{$(1000,1000, 1000, 1000)$}}$,
$g^{(3)}=\mbox{\small{$(1000,1000,4000,4000)$}}$, $g^{(4)}=\mbox{\small{$(4000,4000,1000,1000)$}}$, and $g^{(5)}=\mbox{\small{$(4000,500,4000,500)$}}$.

For the case of adaptive sampling, the initial gain was {\small{$(0,0,0,0)$}}, as before, but since we have more variables and in order to improve robustness, in the generation of the initial subspace we have added a starting point informed by the optimization procedure (that is, at an initial gain of {\small{$(1000,1000, 1000, 1000)$}}).

Note that in these examples the starting point for optimization also included predetermined sampling of damping gain configurations.

 \begin{figure}[h!]
\begin{center}\centering
\includegraphics[width=16cm]{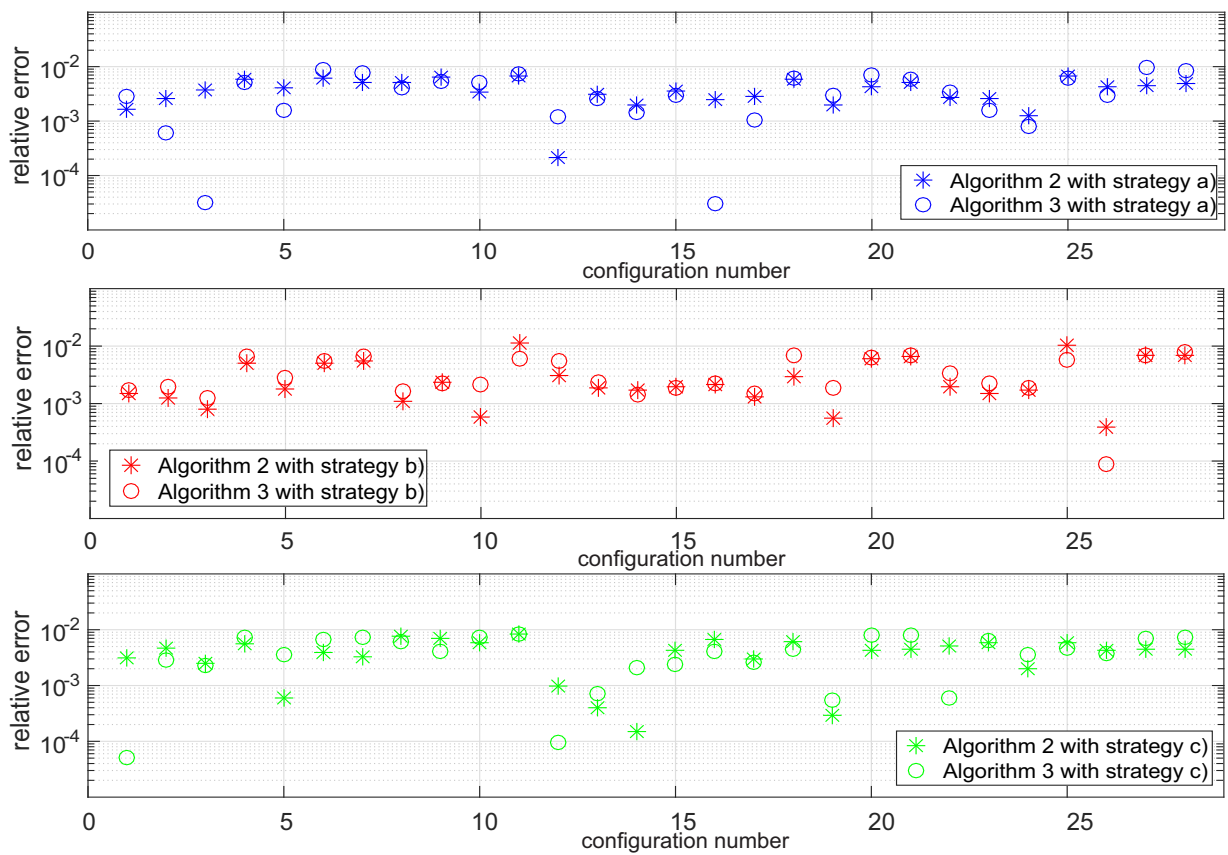}
\caption{Example 2, relative errors for $\cH_2$ norm at optimal gain for Algorithm 2 and Algorithm 3 with strategies a), b) and c)}
 \label{relative-error-abc2001}
\end{center}
\end{figure}

Similarly as in the first example, in Figure \ref{relative-error-abc2001}
 we show relative errors    for the $\cH_2$ norm at optimal gain. This  figure shows results obtained by  Algorithm \ref{OptGain-i} and Algorithm \ref{OptGain-ii} with strategies of  internal reduction based   strategies from cases (a), (b) and (c). We can see that within the given tolerances
  the proposed methodology for damping optimization once again yields satisfactory approximation
  together with significant acceleration of optimization.   As in the previous example,  our approach based on \textsf{sym2IRKA}  returned approximation of optimal gains with better relative errors compared to the one obtained by dominant pole approach based on \cite[Algorithm 1 and 2]{BennerKuerTomljTruh15}  in this example as well.

\begin{table} \begin{center}
\begin{center}
\begin{tabular}{|c||c|c|}
  \hline
     Acceleration factor for  & Algorithm 2 & Algoritm 3 \\ \hline \hline
   Strategy (a)             &     124.72 &  157.98  \\\hline
   Strategy (b)             &     133.53 &  126.16  \\\hline
   Strategy (c)             &     171.93 &  208.80  \\
  \hline
\end{tabular}
\end{center}
\caption{Time ratios for Example 2, i.e., acceleration factor with model reduction accelerated optimization}\end{center}
\end{table}

Similar to the first example, in Table 2 we illustrate average time speed-up for the optimization process obtained using \textsf{sym2IRKA} for different strategies of internal reduction. In this example dominant pole algorithm based on \cite[Algorithm 1 and 2]{BennerKuerTomljTruh15} accelerated  optimization  process 83.87 times. This is  also significant acceleration factor, but new approach based on \textsf{sym2IRKA} was even more efficient, yielding an accelerating factor as high as $208$.  From this table we can see that adaptive strategy is more efficient for strategies a) and b), but  provides slightly smaller time ratio  for strategy b). Here adaptive sampling needed more updates of initial gains which  have had   impact in the final time ratio, but this strategy do not need to have additional information on the area where are  the optimal gains.


Although the time speed ups depend on the tuning tolerances, for both examples we can conclude,  based on performed numerical experiments, that  the approach based on \textsf{sym2IRKA} is significantly faster than the dominant pole approach while also providing  better relative errors.
For the predetermined sampling strategy, we have used gain samples that are somehow in the vicinity of the optimal gains and this is, in general, hard to know.  Thus, in practice, the adaptive sampling strategy is expected to be much more efficient. The adaptive sampling strategy yields similar relative errors as the predetermined sampling strategy and usually provides better or at least similar acceleration speed up.

\section*{Acknowledgements}

 The work of Zoran Tomljanovi\'c was supported in part by the Croatian Science Foundation under the project Optimization of parameter dependent mechanical systems, Grant No. IP-2014-09-9540 and project Control of Dynamical Systems, Grant No. IP-2016-06-2468. The work of Christopher Beattie was supported in part by the Einstein Foundation Berlin.
The work of Serkan Gugercin was supported in part by the Alexander von Humboldt Foundation and by the NSF through Grant DMS-1522616.



\providecommand{\bysame}{\leavevmode\hbox to3em{\hrulefill}\thinspace}
\providecommand{\MR}{\relax\ifhmode\unskip\space\fi MR }
\providecommand{\MRhref}[2]{%
  \href{http://www.ams.org/mathscinet-getitem?mr=#1}{#2}
}
\providecommand{\href}[2]{#2}

\end{document}